# SIR model with random diffusion, reinfection, and random transmission: exponential attractors and spread of the disease

T. Caraballo[a], J. López-de-la-Cruz[b], A. N. Oliveira-Sousa[c],[1]
and P. N. Seminario-Huertas[b]

[a] Dpto. de Ecuaciones Diferenciales y Analysis Numérico, Universidad de Sevilla, Campus Reina Mercedes, Sevilla, Spain.

[a] Department of Mathematics, Wenzhou University, Wenzhou, Zhejiang Province, 325035, P. R. China.

[b] Dpto. de Matemática Aplicada a las TIC, Escuela Técnica Superior de Ingenieros Informáticos, Campus de Montegancedo, Universidad Politécnica de Madrid, 28660 Boadilla del Monte, Madrid, Spain.

[c] Dpto. de Matemática, Universidade Federal de Santa Catarina, 88040-900 Florianópolis, Santa Catarina, Brazil

**Abstract**

We introduce a stochastic SIR-type partial differential equation model incorporating random diffusion, reinfection, vital dynamics, and a randomly varying transmission rate. For the associated random dynamical system, we prove the existence of both random and exponential attractors. We construct a nonstationary, random disease-free global solution, which serves to localize the random attractor. Furthermore, we analyze the mean value of the random transmission coefficient to establish conditions under which the disease may either be eradicated or persist in an endemic state, depending on the system's parameters.

## 1. Introduction

A fundamental model to describe the spread of infectious diseases is the SIR model, introduced by Kermack and McKendrick [15]) in 1927. This model divides the population into three categories: the first category represents susceptible individuals ($S$), who are not currently infected but can contract the disease through contact with an infected individual, the infected individuals ($I$), while the third category comprises recovered individuals ($R$), who either gain immunity after recovering or pass away due to the disease.

This classical model involves the transmission coefficient $\gamma > 0$, which relates $S$ and $I$ by $\gamma SI$, and the mean time $1/c > 0$ of infected $I$ to became recoverd $R$. In recent years, several new parameters have been incorporated into the classical SIR model: a recruitment rate $\Lambda$, which accounts for the influx of new susceptible individuals into the population; a natural death rate $d > 0$, which captures mortality unrelated to the disease; and a reinfection parameter $b > 0$ modeling the loss of immunity in recovered individuals. Several modifications of the classical model were considered with nonautonomous and random/stochasctic models, see [16, 17, 5, 20] and the references therein. In [20], the authors studied the existence of attractors and conditions under which the disease may be eradicated or persist, for this system and several of its variants. These include models in which the transmission coefficient $\gamma$ is treated as a random parameter and the recruitment rate $\Lambda$ is allowed to be time-dependent or stochastic.

---

[1]corresponding author: alexandre.n.o.sousa@ufsc.br

We aim to extend the model introduced in [20] through the incorporation of a random diffusion term, which brings spatial and stochastic dynamics into play, meaning that individuals not only change their state over time but can also move randomly across space. Hence, the populations will also depend on the position in a given bounded domain $\mathcal{O} \subset \mathbb{R}^3$ which we assume has smooth boundary. We model the randomness by parameters on a *metric dynamical system* $(\Omega, \mathcal{F}, \mathbb{P}, \theta)$, i.e., $(\Omega, \mathcal{F}, \mathbb{P})$ is a probability space with a family of measurable maps $\{\theta_t : \Omega \to \Omega; t \in \mathbb{R}\}$ such that $\theta_0 = Id_\Omega$; $\theta_{t+s} = \theta_t \circ \theta_s$, for all $t, s \in \mathbb{R}$ and $\mathbb{P}\theta_t^{-1} = \mathbb{P}$ for all $t \in \mathbb{R}$.

Thus, the resulting modified SIR model is governed by the following system of partial differential equations:

$$S_t = \Lambda + div(a(\theta_t\omega, x)\nabla S) - dS + bI - \gamma(\theta_t\omega)\frac{SI}{S+I+R},\ t \geq t_0,\ x \in \mathcal{O}, \qquad (1)$$

$$I_t = div(a(\theta_t\omega, x)\nabla I) - (d+b+c)I + \gamma(\theta_t\omega)\frac{SI}{S+I+R},\ t \geq t_0,\ x \in \mathcal{O}, \qquad (2)$$

$$R_t = div(a(\theta_t\omega, x)\nabla R) + cI - dR,\ t \geq t_0,\ x \in \mathcal{O}, \qquad (3)$$

with Dirichlet boundary condition, where $S$, $I$, and $R$ denote the size of the population of susceptible, infected, and recovered individuals, respectively, and $\gamma : \Omega \to \mathbb{R}$ and $a(\cdot, x) : \Omega \to \mathbb{R}$ are random variables for all $x \in \mathcal{O}$ and $\Lambda \in \mathbb{R}$.

In [24] they consider a SIR model with a the Laplace operator with Dirichlet Boundary Condition ($a$ and $\gamma$ as constants), and time-dependent recruitment rate, they study the existence of attractors, and the stability of the disease-free (nonautonomous) solution, which leads to conditions under the disease became endemic or not. Recently, in [18] the authors study the existence of random attractors (and finite dimensionality) for stochastic differential equations with the randomness in the linear operator. Motivated by these works, our goal is study the asymptotic behaviour of the model (1)-(3), by means of the theory of random dyamical systems, or more precisely, the theory of random exponential attractors, introduced in [18].

In this work, we investigate the long-term dynamics of (1)–(3) within the framework of random dynamical systems (RDS). Building on the approach of [18], we prove the existence of both a *random attractor* and an *exponential attractor* with finite fractal dimension (see definitions in 2) for the abstract Cauchy problem associated with (1)–(3); see System (4)–(6) in Section 2.

Our analysis further yields conditions that determine whether the disease becomes endemic or whether the infection dies out, based on the asymptotic behavior of a distinguished global solution representing the disease-free state. The theory of random attractors was developed over the past decades and has become a fundamental tool for studying infinite-dimensional dynamical systems. The pioneering works are [11, 12], were followed by a range of studies addressing the existence, stability, and structure of random attractors, such as [7, 8, 9, 26]. For a comprehensive introduction to this topic, we refer to [6]. Subsequently, the concept of an exponential random attractor was introduced in [10], which relaxes the negative invariance property but ensures an exponential rate of attraction and finite fractal dimension. Dimension theory for attractors is crucial because embedding results allow infinite-dimensional dynamics to be represented in finite-dimensional spaces, facilitating both theoretical analysis and numerical approximation (see, e.g., [13, 21, 23]). In the setting of random dynamical systems, the abstract framework for random exponential attractors is developed in [10], with some applications discussed in [18].

We show that the abstract SIR model generates a random dynamical system admitting an exponential random attractor. Moreover, we propose a way to study the spread of the disease using the mean-time value in time of the stochastic process $\gamma(\theta_t\omega)$, see (14). Then we investigate conditions under which the disease is eradicated ( the number of infected tends to zero, see Theorem 4.1) or becomes endemic (the infected population persists in a weak sense, see Theorem 4.2). This analysis essentially reduces to studying the stability of a non-stationary disease-free solution: if the solution is exponentially stable, the attractor coincides with it, and the infected population vanishes asymptotically; if the solution is unstable, the disease persists endemically. In the latter case, it is of particular interest to examine the fractal dimension and structure of the resulting random attractor.

## 2. Preliminaries and Problem Formulation

We begin by reformulating the random SIR model introduced in the previous section within the abstract framework of random dynamical systems. This reformulation allows us to treat the problem from the viewpoint of infinite-dimensional dynamics with coefficients driven by an underlying metric dynamical system, and it provides the foundation for the global analysis developed later in the paper.

The system (1)–(3) naturally lead to the formulation of the SIR model as an abstract

Cauchy problem

$$\dot{S} = \Lambda + A(\theta_t\omega)S - dS + bI - \gamma(\theta_t\omega)\frac{SI}{S+I+R}, \quad (4)$$

$$\dot{I} = A(\theta_t\omega)I - (d+b+c)I + \gamma(\theta_t\omega)\frac{SI}{S+I+R}, \quad (5)$$

$$\dot{R} = A(\theta_t\omega)R + cI - dR, \quad (6)$$

in $X_+ = (L^2_+(\mathcal{O}))^3$, where $X = L^2(\mathcal{O})^3$ and $L^2_+(\mathcal{O})$ the subspace of non-negative functions of $L^2(\mathcal{O})$ and $\|(u_1, u_2, u_3)\|_X = \|u_1\|_2 + \|u_2\|_2 + \|u_3\|_2$.

**2.1. Analytical assumptions and existence of solutions**

In line with the framework presented in [18], we consider $\{A(\omega) : \omega \in \Omega\}$ to be a family of operators satisfying the following:

(A0) the operators $A(\omega) : D \subset L^2(\mathcal{O}) \to L^2(\mathcal{O})$ are closed, densely defined and have a common domain $D = H^2(\mathcal{O}) \cap H^1_0(\mathcal{O})$, for all $\omega \in \Omega$;

(A1) the mapping $(t, \omega) \mapsto A(\theta_t\omega) \in \mathcal{L}(D, L^2(\mathcal{O}))$ is strongly measurable;

(A2) there exists $\vartheta \in (\pi, \pi/2)$ and $M > 0$ such that $\Sigma_\vartheta = \{\mu \in \mathbb{C} : |arg(\mu)| \leq \vartheta\} \subset \rho(A(\omega))$ and

$$\|(\mu - A(\omega))^{-1}\|_{\mathcal{L}(L^2(\mathcal{O}))} \leq \frac{M}{|\mu|+1}, \ \mu \in \Sigma_\vartheta \cup \{0\}, \ \omega \in \Omega. \quad (7)$$

(A3) There exists $\nu \in (0, 1]$ and $C > 0$ such that

$$\|A(\theta_t\omega) - A(\theta_s\omega)\|_{\mathcal{L}(D, L^2(\mathcal{O}))} \leq C|t-s|^\nu \ \ \forall \ t, s \in \mathbb{R}, \ \omega \in \Omega; \quad (8)$$

(A4) The operators $A(\omega)$ have compact inverse for all $\omega \in \Omega$;

(A5) there exists $\lambda_0 > 0$ such that $(-A(\omega)v, v)_2 \geq \lambda_0 \|v\|_2^2$, for all $v \in D$.

Note that, conditions (A2)–(A3) are the so-called *Kato–Tanabe conditions*, which guarantee the existence of an evolution family associated with the nonautonomous operator $A(\theta_t\omega)$, see for instance [2, p. 55]. The sectoriality and resolvent bound in (A2) ensure that each $A(\omega)$ generates an analytic semigroup uniformly in $\omega$, providing parabolic-type smoothing properties even when coefficients are random, on the other hand the condition (A3) — a Hölder-type regularity in time — is the key assumption that allows the construction of an evolution family $U(t, s, \omega)$ via the Acquistapace–Terreni theory and Kato-Tanabe theory, since the system is nonautonomous. Compactness (A4) yields precompact orbits, crucial for attractor theory, while coercivity (A5) provides uniform exponential dissipation. Examples satisfying (A0)–(A5) include

$$A_1(\omega) = a(\omega)\Delta, \quad A_2(\omega) = a(\omega) + \Delta, \quad A(\omega) = \text{div}(a(\omega, x)\nabla \cdot),$$

as proved in [18]. The diffusion operator appearing in (1)–(3) is also of this type, and its specific role in the dynamics of disease spread will be further analyzed in Section 4.

In [18, Theorem 2.14], it is explained how to use [1, Theorem 2.3] and [22, Theorem 2.2] to obtain the following result.

**Theorem 2.1** (Random Hille–Yosida type result)**.** *Assume* (A0)–(A5)*. Then for every* $\omega \in \Omega$ *there exists an evolution family* $U(t,s,\omega) : L^2(\mathcal{O}) \to L^2(\mathcal{O})$, $t \geq s$, *such that:*

1. $U(t,t,\omega) = \mathrm{Id}$ *and* $U(t,s,\omega)U(s,\tau,\omega) = U(t,\tau,\omega)$ *for all* $t \geq s \geq \tau$;

2. $U(t,s,\omega)$ *is strongly continuous and satisfies*

$$\partial_t U(t,s,\omega) = A(\theta_t \omega) U(t,s,\omega) \qquad \text{for } t > s;$$

3. *There exists* $\lambda_0 > 0$ *such that*

$$\|U(t,s,\omega)\|_{\mathcal{L}(L^2(\mathcal{O}))} \leq e^{-\lambda_0 (t-s)}, \qquad t \geq s.$$

Now, we use this evolution family to study nonautonomous random problems of the form

$$\dot{u}(t) = A(\theta_t \omega) u(t) + F(\theta_t \omega, u(t)), \qquad u(t_0) = u_0 \in X, \tag{9}$$

where $F : \Omega \times X \to X$ is the nonlinearity.

The next lemma, adapted from Lemma 3.6 in [18], summarizes the corresponding well-posedness theory.

**Lemma 2.1** (Well-posedness for random nonautonomous evolution equations)**.** *Let* (A0)–(A5) *hold, and let* $F : \Omega \times X \to X$ *satisfy:*

(i) *for each* $\omega \in \Omega$, *the map* $u \mapsto F(\omega, u)$ *is globally Lipschitz on* $X$, *with constant* $L_F$ *independent of* $\omega$;

(ii) *for each* $u \in X$, *the map* $\omega \mapsto F(\omega, u)$ *is measurable.*

*Then for every* $t_0 \in \mathbb{R}$, $\omega \in \Omega$ *and* $u_0 \in X$, *there exists a unique mild solution*

$$u(\cdot, t_0, \omega; u_0) \in C([t_0, +\infty); X)$$

*of* (9)*, given by the variation-of-constants formula*

$$u(t, t_0, \omega; u_0) = U(t, t_0, \omega) u_0 + \int_{t_0}^{t} U(t, s, \omega) F(\theta_s \omega, u(s, t_0, \omega; u_0)) \, ds.$$

*Moreover, the mapping*

$$\psi(t, t_0, \omega) u_0 := u(t, t_0, \omega; u_0)$$

*defines a continuous evolution process on* $X$.

Now, we use Theorem 2.1 to define mild solutions for (4)-(6). Let $t_0 \in \mathbb{R}$, $u_0 \in X_+$, and $\omega \in \Omega$ a **mild solution** for (4)-(6) through $u_0$ at $t_0$ is a continuous function $[t_0, \sigma] \ni t \mapsto u(t, t_0, \omega, u_0) \in D$, given by

$$u(t, t_0, \omega, u_0) = (S(t, t_0, \omega, u_0), I(t, t_0, \omega, u_0), R(t, t_0, \omega, u_0)),$$

where $S, R$ and $I$ satisfy

$$\begin{aligned}
S(t, t_0) &= U(t, t_0, \omega)S_0 + \int_{t_0}^{t} U(t, r, \omega)(\Lambda - dS(r, t_0) + bI(r, t_0))\, dr \\
&\quad - \int_{t_0}^{t} U(t, r, \omega)\left(\gamma(\theta_r \omega)\frac{S(r, t_0)I(r, t_0)}{N(r, t_0)}\right) dr, \\
I(t, t_0) &= U(t, t_0, \omega)I_0 + \int_{t_0}^{t} U(t, r, \omega)\left(-(d+b+c)I(r, t_0) + \gamma(\theta_r \omega)\frac{S(r, t_0)I(r, t_0)}{N(r, t_0)}\right) dr, \\
R(t, t_0) &= U(t, t_0, \omega)R_0 + \int_{t_0}^{t} U(t, r, \omega)(cI(r, t_0) - dR(r, t_0))\, dr.
\end{aligned} \tag{10}$$

Above, we abused notation by writing $(t, t_0)$ to mean $(t, t_0, \omega, u_0)$. This simplification will be used throughout when no confusion arises. If the maximal time of existence is $\sigma = +\infty$, for each $\omega \in \Omega$ fixed, the solution $u$ is associated with the following **evolution process**

$$\psi(t, t_0, \omega)u_0 := u(t, t_0, \omega, u_0),$$

i.e., satisfies $\psi(\cdot, \cdot, \omega)$ and $\psi(t, t, \omega) = Id_X$, for all $t \in \mathbb{R}$; $\psi(t, s, \omega)\psi(s, \tau, \omega) = \psi(t, \tau, \omega)$, for all $t \geq s \geq \tau$ and $\omega \in \Omega$; and the mapping $(t, t, t_0, u_0) \mapsto \psi(t, t_0, \omega)u_0$ is strongly continuous.

Let $\{W_t\}_{t \in \mathbb{R}}$ be the standard Wiener process, i.e., a stochastic process considered defined in a probability space $(\Omega, \mathcal{F}, \mathbb{P})$ is defined as follows: $\Omega$ is the space of continuous functions vanishing at $t = 0$, denoted by $C_0(\mathbb{R})$, with $W_t(\omega) = \omega(t)$ for all $\omega \in C_0(\mathbb{R})$ and $t \in \mathbb{R}$. The space $\Omega$ is equipped with the compact-open topology, and $\mathcal{F}$ represents the Borel $\sigma$-algebra generated by this topology. The probability measure $\mathbb{P}$ on $\mathcal{B}(C_0(\mathbb{R}))$ is the Wiener measure corresponding to the distribution of a two-sided Wiener process with trajectories in $C_0(\mathbb{R})$. Using the Wiener shift, defined as $\theta_t \omega := \omega(t + \cdot) - \omega(t)$ for $t \in \mathbb{R}$ and $\omega \in \Omega$, the resulting structure $(\Omega, \mathcal{F}, \mathbb{P}, (\theta_t)_{t \in \mathbb{R}})$ provides a metric dynamical system, see [6].

In [18, Theorem 3.1], it is shown that if $\Omega$ is the canonical sample space of the Wierner process, then $(t, \omega) \mapsto U(t, 0, \omega)$ defines a random dynamical system driven by the Wierner shift $(\theta_t)_{t \in \mathbb{R}}$. Thus problem (10) is associated with a **random dynamical system (RDS)** in $X$:

$$\varphi(t, \omega)u_0 := u(t, 0, \omega; u_0),\ t \geq 0, \omega \in \Omega \text{ and } u_0 \in X,$$

i.e., $\varphi(0, \omega)u_0 = u_0$ for every $u_0 \in X$ and $\omega \in \Omega$, $\varphi(t + s, \omega) = \varphi(t, \theta_s \omega)\varphi(s, \omega)$ for $t \geq s$; $\omega \mapsto \varphi(t, \omega)u_0$ is measurable for all $t, u_0$ fixed; and $\varphi(t, \omega) : X \to X$ is continuous.

Hence, if $\Omega$, since $u(t + t_0, t_0, \omega; u_0) = u(t, 0, \theta_{t_0}\omega, u_0)$, we have the following relation between the evolution process $\psi$ and the random dynamical system $\varphi$:

$$\psi(t, t_0, \omega) = \varphi(t - t_0, \theta_{t_0}\omega), \ t \geq t_0, \ \omega \in \Omega.$$

We now, briefly recall the concepts of attractors for RDS.

A **random attractor** for a RDS $\varphi$ is a family $\{\mathcal{A}(\omega) : \omega \in \Omega\}$ of compact subsets of $X_+$ such that:

1. $\{\mathcal{A}(\omega) : \omega \in \Omega\}$ is a random set, i.e., $x \mapsto d(x, \mathcal{A}(\omega))$ is measurable for all $x \in X_+$.

2. $\varphi(t, \omega)\mathcal{A}(\omega) = \mathcal{A}(\theta_t\omega), \ t \geq 0$.

3. $\mathcal{A}$ is pullback attracts bounded subsets, i.e., given a bounded subset $D$ of $X_+$ we have
$$\lim_{t \to +\infty} dist(\varphi(t, \theta_{-t}\omega)D, \mathcal{A}(\omega)) = 0,$$
where $dist(A, B) = \sup_{a \in A} \inf_{b \in B} \|a - b\|_X$ denotes the Hausdorff semi-distance between $A, B \subset X_+$.

A **random exponential attractor** is a family of compact sets $\{\mathcal{M}(\omega) : \omega \in \Omega\}$ such that:

1. $\{\mathcal{M}(\omega) : \omega \in \Omega\}$ is a random set;

2. $\mathcal{M}$ is positively invariant, i.e., $\varphi(t, \omega)\mathcal{M}(\omega) \subset \mathcal{M}(\theta_t\omega), \ t \geq 0, \ \omega \in \Omega$;

3. $\mathcal{M}$ is pullback attracting with an exponential rate, i.e., there exists $\alpha > 0$ such that for any bounded subset $B$ of $X_+$,
$$\lim_{t \to +\infty} e^{\alpha t} dist(\varphi(t, \theta_{-t}\omega)B, \mathcal{M}(\omega)) = 0;$$

4. $\mathcal{M}$ has finite fractal dimension, i.e., there exists $k > 0$ such that $dim_F(\mathcal{M}(\omega)) \leq k$, for all $\omega \in \Omega$,

where the **fractal dimension** of a precompact subset $K$ of $X$ is defined as

$$dim_F(K) = \limsup_{\epsilon \to 0} \log_{1/\epsilon}(N_\epsilon(K)),$$

where $N_\epsilon(K)$ denotes the minimal number of balls of radius $\epsilon > 0$ in $K$ needed to cover the set $K$.

### 3. Existence of random exponential attractors

We study the existence and uniqueness of positive solutions of (4)-(6) and continuity with respect to the initial data in $X_+ = (L_+^2(\mathcal{O}))^3$, where $L_+^2(\mathcal{O})$ the subspace of non-negative functions of $L^2(\mathcal{O})$.

Caraballo et. al.**Proposition 3.1.** *For any $u_0 = (S_0, I_0, R_0) \in X_+$ and $\omega \in \Omega$, there exists a unique solution $u(\cdot, t_0, \omega; u_0) : [0, +\infty) \to X_+$ of (4)-(6), which is associated to a RDS $\varphi : [t_0, +\infty) \times \Omega \times X_+ \to X_+$ defined as $\varphi(t, \omega)u_0 = u(t + t_0, t_0, \omega; u_0)$ for all $t \geq 0$, $u_0 \in X_+$ and $\omega \in \Omega$. In addition, there exists a bounded set subset of $X_+$, given by*

$$B_1 = \left\{ (S, I, R) \in X_+ : \|S + I + R\|_2 \leq 1 + \frac{\|\Lambda\|_2}{\lambda_0 + d} \right\},$$

*such that absorbs bounded subsets of $X_+$, in the pullback and forward sense.*

*Proof.* The proof of well-posedness, and that the solutions generates a RDS in $X_+$ are standard, following similar steps of [18, Lemma 3.6 ] and [14, Theorem 2.1].

Let $t_0 \in \mathbb{R}$ and $\omega \in \Omega$ be fixed. Let us analyze the equation for $N = S + I + R$, we have

$$N_t = \Lambda + A(\theta_t \omega) N - dN, \ t \geq t_0, \ N(t_0) = N_0. \tag{11}$$

From Theorem 2.1 and [19, Lemma 5.3], for $t \geq t_0$, $N$ is given by

$$N(t, t_0, \omega, N_0) = U(t, t_0, \omega) e^{-d(t-s)} N_0 + \int_{t_0}^{t} U(t, r, \omega) e^{-d(t-r)} \Lambda \, dr, \text{ for }, N_0 \in L_+^2(\mathcal{O}),$$

Thus,

$$\begin{aligned} \|N(t, t_0, \omega; N_0)\|_2 &\leq e^{-(\lambda_0 + d)(t - t_0)} \|N_0\|_2 + \int_{t_0}^{t} e^{-(\lambda_0 + d)(t-r)} \|\Lambda\|_2 \, dr \\ &\leq e^{-(\lambda_0 + d)(t - t_0)} \|N_0\|_2 + \frac{\|\Lambda\|_2}{\lambda_0 + d}. \end{aligned} \tag{12}$$

From above, we obtain that the solution $u(t, t_0, \omega; u_0)$ is defined for all $t \geq t_0$, and that $B_1$, defined in (3.1), is a bounded subset of $X_+$ that pullback absorbs bounded subsets. $\square$

To apply the results of [18] we study the nonlinearity from (4)-(6).

**Lemma 3.1.** *Define $f : L_+^2(\mathcal{O}) \setminus \{0\} \times L^2(\mathcal{O}) \times L^2(\mathcal{O}) \to L^2(\mathcal{O})$ by*

$$f(S, I, R) = \frac{SI}{S + I + R}.$$

*Then $f$ is Fréchet differentiable with bounded derivative. In particular the mapping*

$$(S, I, R) \mapsto (\gamma f(S, I, R), -\gamma f(S, I, R), 0)$$

*is globally Lipschitz with Lipschitz constant $C_f(\gamma) \leq 2\gamma + 1$, for $\gamma > 0$.*

*SIS-IR model with random diffusion and vital dynamics*

*Proof.* It is straightforward to verify that $f$ is well defined and that the partial derivatives of $f$, given by

$$\partial_S f(S, I, R) = \frac{I(R+I)}{N^2},$$
$$\partial_I f(S, I, R) = \frac{S(R+S)}{N^2},$$
$$\partial_R f(S, I, R) = \frac{-SI}{N^2},$$

are well-defined linear operators, where $N = S + I + R$. Thus the gradient $\nabla f : L^2_+(\mathcal{O}) \setminus \{0\} \times L^2(\mathcal{O}) \times L^2(\mathcal{O}) \to L^2(\mathcal{O})$ is such that

$$\|\nabla f(S, I, R)\|_{\mathcal{L}(L^2_+(\mathcal{O}) \times L^2(\mathcal{O}) \times L^2(\mathcal{O}))} \leq 3, \quad \text{for all } (S, I, R).$$

The last claim follows from the mean-value inequality, and the proof is complete. $\square$

As a consequence of Proposition 3.1, Lemma 3.1 and Theorem 3.8 from [18], in the following result we obtain the existence of exponential attractor for (4)-(6).

**Theorem 3.1.** *Assume that $\gamma : \Omega \to [0, +\infty)$ is bounded. Then for every $\nu \in (0, 1/2)$ and $\eta \in (0, 1)$, there exists a random pullback exponential attractor $\mathcal{M} := \mathcal{M}^{\nu,\eta}$ satisfying,*

$$dim_F(\mathcal{M}(\omega)) \leq \log_{\frac{1}{2\nu}} \left( N^X_{\frac{\kappa}{2\nu}}(B^{X^\eta}(0,1)) \right), \ \nu \in (0, 1/2).$$

*Moreover, there exists a unique random attractor and its fractional dimension is bounded by*

$$dim_F(\mathcal{A}(\omega)) \leq \inf_{\nu \in (0, \frac{1}{2})} \left\{ \log_{\frac{1}{2\nu}} \left( N^X_{\frac{\kappa}{2\nu}}(B^{X^\eta}(0,1)) \right) \right\}.$$

*Proof.* We have to verify the smoothing property needed in the abstract results on the existence of exponential attractors, see Conditions (H3) and (H4) of [18, Thm 2.2]. Indeed, note that from Lemma 3.1 together with the condition that $\gamma$ is bounded, it follows that the nonlinearity of globally Lipschitz and then we apply Theorem 3.8 from [18]. $\square$

**Remark 3.1.** *The result [18, Thm 2.2] is a simpler case of the main result of [10]. In [10] it is possible to obtain $\kappa$ above as a random variable. In our case is simpler because of the property of the nonlinearity of (4)-(6).*

**Definition 3.1.** *We say that $\{\mathcal{A}(\omega) : \omega \in \Omega\}$ is a **bounded random attractor** if*

$$\bigcup_{\omega \in \Omega} \mathcal{A}(\omega) \text{ is a bounded subset of } X.$$

**Remark 3.2.** *It is straightforward to verify that if $\{\mathcal{A}(\omega) : \omega \in \Omega\}$ is bounded, then*

$$\mathcal{A}(\omega) = \bigcup \{\xi(\omega) : \xi : \Omega \to X \text{ is a bounded global solution of } \varphi\}, \tag{13}$$

*Caraballo et. al.*

*see [4, Theorem 6.11] for the nonautonomous case.*

*If $\Lambda$ is bounded the random attractor of (10) is bounded, so it can be characterized as (13).*

From the equation for $N = S + I + R$, it appears a special solution.

**Theorem 3.2.** *Define*

$$N^*(t,\omega) := \int_{-\infty}^{t} U(t,r,\omega)e^{-d(t-r)}\Lambda\, dr, \ (t,\omega) \in \mathbb{R} \times \Omega,$$

*with this limit being taken in $L^2(\mathcal{O})$. Then $N^*$ defines a bounded global solution for (11) and*

$$\Sigma^*(t,\omega) := \{(S,I,R) \in X_+ : S + I + R = N^*(t,\omega)\}$$

*is a bounded subset that attracts pullback and forward bounded subsets of $X_+$. In particular, $\mathcal{A}(\omega) \subset \Sigma^*(0,\omega)$, for all $\omega \in \Omega$.*

*Proof.* In fact, note that

$$N^*(t,\omega) = \lim_{s\to-\infty} N(t,s,\omega,N_0), \text{ for all } N_0 \in L^p(\mathcal{O}).$$

Now, from (12), we see that $N^*$ is a bounded global solution, $\|N^*(t,\omega)\| \leq \frac{\|\Lambda\|_p}{\lambda_0+d}$. Also, since $N(t,s,\omega,N^*(s,\omega)) = N^*(t,\omega)$, for $t > s$, we have

$$N^*(t,\omega) = U(t,s,\omega)e^{-d(t-s)}N^*(s,\omega) + \int_s^t U(t,r,\omega)e^{-d(t-r)}\Lambda\, dr, \ t \geq s.$$

Hence,
$$N(t,s,\omega,N_0) - N^*(t,\omega) = U(t,s,\omega)e^{-d(t-s)}(N_0 - N^*(s,\omega)),$$

and from (A4), this goes to zero in $L^p(\mathcal{O})$ as $t \to +\infty$, for each $s$ and $\omega$ fixed. For the last claim use the fact that the random attractor is bounded and invariant, and the proof is complete. □

**Remark 3.3.** *One can prove the existence of attracts directly, see [24]. The candidate will be the pullback limit set of the absorbing set, $B_1$ given in (3.1), i.e.,*

$$K(\omega) := \bigcap_{t \geq 0} \overline{\varphi(t,\theta_{-t}\omega)B_1}$$

**4. Spread of the disease**

In this section, we propose a way to study the spread of the disease using the mean-time value of the stochastic process $\gamma(\theta_t\omega)$.

Define
$$m(\gamma,\omega) := \lim_{n\to+\infty} \sup\left\{\frac{1}{t-t_0}\int_{t_0}^t \gamma(\theta_r\omega)\,dr : t-t_0 > n\right\}, \quad \text{for all } \omega \in \Omega. \tag{14}$$

Now, we study a condition for the disease to be eradicated. The following result is inspired on the finite-dimensional case, see [20].

**Theorem 4.1.** *Assume that*
$$m(\gamma,\omega) < \lambda_0 + d + b + c, \quad \text{for all } \omega \in \Omega. \tag{15}$$

*Then, for all $\omega \in \Omega$,*

$$\|(S(t,t_0,\omega,u_0), I(t,t_0,\omega,u_0), R(t,t_0,\omega,u_0)) - (N^*(t,\omega),0,0)\|_X \to 0, \quad \text{as } t-t_0 \to +\infty.$$

*In particular, the random attractor is given by $\mathcal{A}(\omega) = \{(N^*(0,\omega),0,0)\}$.*

*Proof.* Let $u_0 = (S_0, I_0, R_0) \in X_+$, $\omega \in \Omega$, and $\hat{U}(t,t_0) := U(t,t_0,\omega)e^{-(d+b+c)(t-t_0)}$ and note that $I$ satisfies

$$I(t,t_0,\omega;u_0) = \hat{U}(t,t_0)I_0 + \int_{t_0}^t \hat{U}(t,r)\gamma(\theta_r\omega)\frac{S(r,t_0,\omega;u_0)I(r,t_0,\omega;u_0)}{N(r,t_0,\omega;u_0)}\,dr, \tag{16}$$

see for instance [19, Lemma 5.3]. Since $\|\hat{U}(t,t_0)\|_{\mathcal{L}(X)} \leq e^{-(\lambda_0+d+b+c)(t-t_0)}$ and $\|SI/N\|_2 \leq \|I\|_2$, we have

$$\|I(t,t_0,\omega;u_0)\|_2 \leq e^{-(\lambda_0+d+b+c)(t-t_0)}\|I_0\|_2 + \int_{t_0}^t e^{-(\lambda_0+d+b+c)(t-r)}\gamma(\theta_r\omega)\|I(r,t_0,\omega;u_0)\|_2\,dr.$$

Applying Grönwall's inequality for $\varphi(t) = \|I(t,t_0,\omega;u_0)\|_X e^{(\lambda_0+d+b+c)t}$, we obtain

$$\|I(t,t_0,\omega;u_0)\|_2 \leq e^{-(\lambda_0+d+b+c)(t-t_0)+\int_{t_0}^t \gamma(\theta_r\omega)\,dr}\|I(t_0,t_0,\omega;u_0)\|_2. \tag{17}$$

Thus, $\lim_{t-t_0\to+\infty} I(t,t_0,\omega;u_0) = 0$. Now, similarly $R$ satisfies

$$R(t,t_0,\omega;u_0) = U(t,t_0,\omega)e^{-d(t-t_0)}R_0 + \int_{t_0}^t U(t,r,\omega)e^{-d(t-r)}cI(r,t_0,\omega;u_0)\,dr. \tag{18}$$

Now, let us prove that

$$\lim_{t\to+\infty}\int_{t_0}^t U(t,r,\omega)e^{-d(t-r)}cI(r,t_0,\omega;u_0)\,dr = 0,$$

for each $t_0 \in \mathbb{R}$ and $\omega$ fixed.

Indeed, using (17) and that $\|U(t, t_0, \omega)\|_2 \leq e^{-\lambda_0(t-t_0)}$, then

$$\| \int_{t_0}^{t} U(t, r, \omega) e^{-d(t-r)} cI(r, t_0, \omega; u_0) \, dr \|_2$$
$$\leq c \int_{t_0}^{t} e^{-(\lambda_0+d)(t-r) - (\lambda_0+d+c+b)(r-t_0) + \int_{t_0}^{r} \gamma(\theta_s \omega) \, ds} \, dr \quad (19)$$
$$= c e^{-(\lambda_0+d)(t-t_0)} \int_{t_0}^{t} e^{-(c+b)(r-t_0) + \int_{t_0}^{r} \gamma(\theta_s \omega) \, ds} \, dr.$$

For each $t_0$ and $\omega$ fixed, the integral in the right-hand side either exists or goes to infinite. If it does exist, then the limit goes to zero because of exponential decay. If it is infinite, we can apply L'Höpital's rule and obtain

$$\frac{1}{\lambda + d} e^{-(\lambda_0+d)(t-t_0)} e^{-(c+b)(t-t_0) + \int_{t_0}^{t} \gamma(\theta_s \omega)} \, ds \to 0 \text{ as } t \to +\infty,$$

by Condition (15), and hence, either case we obtain $R(t, t_0, \omega; u_0) \to 0$ as $t \to +\infty$. Now we prove that $R(t, t_0, \omega; u_0) \to 0$ as $t_0 \to -\infty$. By the same reasoning above, it is enough to show that

$$\lim_{t_0 \to -\infty} \int_{t_0}^{t} U(t, r, \omega) e^{-d(t-r)} cI(r, t_0, \omega; u_0) \, dr = 0.$$

From Condition (15), there exist $n_0 \in \mathbb{N}$, such that

$$\int_{t_0}^{r} \gamma(\theta_s \omega) \, ds < m(\gamma, \omega)(r - t_0), \; t_0 < -M.$$

Hence from (19), for $t_0 < -M$, we have

$$\| \int_{t_0}^{t} U(t, r, \omega) e^{-d(t-r)} cI(r, t_0, \omega; u_0) \, dr \|$$
$$\leq c e^{-(\lambda_0+d)(t-t_0)} \int_{t_0}^{t} e^{-(c+b)(r-t_0) + \int_{t_0}^{r} \gamma(\theta_s \omega) \, ds} \, dr$$
$$\leq c e^{-(\lambda_0+d)(t-t_0)} \int_{t_0}^{t} e^{-(c+b-m(\gamma,\omega))(r-t_0)} \, dr \to 0, \text{ as } t_0 \to -\infty,$$
$$(20)$$

where the last limit holds true, since $m(\gamma, \omega) < \lambda_0 + d + c + b$.

The remainder of the proof follows from Theorem 3.2.

$\square$

Now, we study a condition for weak uniform persistence. For constant $\gamma$, with $\gamma > \lambda_1 + a + b + c$ and $\Lambda$ is time-dependent and bounded in [3, Thm. 5.1] it is proven the "weak uniform persistence", i.e., for all non-negative initial data $u_0 = (S_0, I_0, R_0)$ with $I_0 > 0$.

For the final result, we will focus on the system (1)-(3), the conditions on $a$ are described bellow.

**Example 4.1.** *Let $a : \Omega \times \mathcal{O} \to \mathbb{R}$ be a random variable such that:*

1. *there are $0 < a_0 < a_1$ such that $a(\omega, x) \in (a_0, a_1)$, for all $i = 1, 2, 3$, $x \in \mathcal{O}$, $\omega \in \Omega$;*

2. *the mapping $t \mapsto \partial_x a(\theta_t \omega, x)$ is continuous for each $\omega \in \Omega$ and $t \in \mathbb{R}$;*

3. *there exists $c > 0$ such that*

$$\sum_{j=1}^{3} a(\omega, x) \xi_j^2 \geq c |\xi|^2,$$

4. *$a_j(\omega, \cdot) \in C^3(\bar{\mathcal{O}})$ for all $\omega \in \Omega$ and there exists $\nu \in (0, 1]$ and $c_1 > 0$, such that*

$$|a_j(\theta_t \omega, x) - a_j(\theta_s \omega, x)| \leq c_1 |t - s|^\nu;$$

*for all $t, s \in \mathbb{R}$, $x \in \overline{\mathcal{O}}$, $\omega \in \Omega$.*

In [18] is shown that all the following operators $A_1(\omega) := a(\omega)\Delta$, $A_2(\omega) := a(\omega) + \Delta$, and $A(\omega) := div(a(\omega, x)\nabla \cdot)$ fulfills (A0)-(A4). Note that $A$ satisfies (A5). Indeed, by Poincaré inequality,

$$\lambda_1 \|u\|_2^2 \leq \|\nabla u\|_2^2, \ u \in H_0^1(\mathcal{O}),$$

where $\lambda_1 > 0$ is the first eigenvalue of $-\Delta$ with Dirichlet boundary condition. Hence, for $u \in H_0^1(\mathcal{O})$, we have

$$(-A(\omega)u, u)_2 = \int_{\mathcal{O}} a(\omega) |\nabla u|^2 dx \geq a_0 \|\nabla u\|_2^2 \geq a_0 \lambda_1 \|u\|_2^2.$$

**Theorem 4.2.** *Let $A(\omega) := div(a(\omega, x)\nabla \cdot)$, where $a : \Omega \times \mathcal{O} \to (a_0, a_1)$ satisfies the conditions on Example 4.1 and $\lambda_1$ is the first eigenvalue of $-\Delta$. Assume that $\gamma$ satisfies*

$$m(\gamma, \omega) > \lambda_1 a_1 + d + b + c, \ \ \text{for all } \omega \in \Omega. \tag{21}$$

*Then,*

$$\limsup_{t \to +\infty} \int_{\mathcal{O}} I(t, t_0, \omega, x) \, dx > 0, \ \forall t_0 \in \mathbb{R} \text{ and } \omega \in \Omega. \tag{22}$$

*Proof.* Let $\omega \in \Omega$ be fixed. Assume that there exists $u_0 = (S_0, I_0, R_0)$ with $I_0 > 0$ and that the solution $u(t, t_0, \omega; u_0) = (S(t, t_0, \omega, x), I(t, t_0, \omega, x), R(t, t_0, \omega, x))$ is such that

$$\limsup_{t \to +\infty} \int_{\mathcal{O}} I(t, t_0, \omega, x) \, dx = 0. \tag{23}$$

$\alpha = d + b + c$, by means of Green's function and applying [25, Lemma 6.4], we are able to prove that given $\epsilon > 0$ with

$$\epsilon < \frac{m(\gamma, \omega) - \lambda_1 - \alpha}{m(\gamma, \omega)}, \tag{24}$$

there exists $\tau_0 = \tau_0(\epsilon, t_0, \omega) \geq t_0$ and such that

$$\frac{I(t, t_0, x) + R(t, t_0, x)}{N(t, t_0, x)} \leq \epsilon, \text{ for } t \geq \tau_0 \text{ and } x \in \mathcal{O}. \tag{25}$$

Let $v_1$ be an eigenvector associated with $\lambda_1$ with $\|v_1\|_2 = 1$. Define

$$w(t, t_0, \omega) := \int_{\mathcal{O}} I(t, t_0, \omega, x) v_1(x) \, dx, \ t \geq t_0.$$

Since $I$ satisfies (5), we have

$$\dot{w}(t, t_0, \omega) = \int_{\mathcal{O}} A(\theta_t \omega) I(t, t_0, \omega, x) v_1 \, dx - \alpha w(t) \tag{26}$$

$$+ \gamma(\theta_t \omega) \int_{\mathcal{O}} \frac{S(t, t_0, \omega, x) I(t, t_0, \omega, x) v_1(x)}{N(t, t_0, \omega, x)} \, dx. \tag{27}$$

First, integrating by parts and using the boundary conditions, we have

$$\begin{aligned}
\int_{\mathcal{O}} A(\theta_t \omega) I(t, t_0, \omega, x) v_1 \, dx &= -\int_{\mathcal{O}} a(\theta_t \omega, x)(\nabla I(t, t_0, \omega, x)) \cdot (\nabla v_1(x)) \, dx \\
&\geq -a_1 \int_{\mathcal{O}} (\nabla I(t, t_0, \omega, x)) \cdot (\nabla v_1(x)) \, dx \\
&= a_1 \int_{\mathcal{O}} I(t, t_0, \omega, x) \Delta v_1(x) \, dx \\
&\geq -a_1 \lambda_1 w(t, t_0, \omega).
\end{aligned} \tag{28}$$

Hence, from (26) we have

$$\dot{w}(t, t_0, \omega) \geq -\lambda_1 a_1 w(t) - \alpha w(t) + \gamma(\theta_t \omega) \int_{\mathcal{O}} \frac{S(t, t_0, x) I(t, t_0, x) v_1(x)}{N(t, t_0, x)} \, dx.$$

Using that $S = N - I - R$, we have

$$\begin{aligned}
\dot{w}(t, t_0, \omega) = &(-a_1 \lambda_1 - \alpha + \gamma(\theta_t \omega)) w(t, t_0, \omega) \\
&- \gamma(\theta_t \omega) \int_{\mathcal{O}} \frac{(I(t, t_0, \omega, x) + R(t, t_0, \omega, x)) I(t, t_0, \omega, x) v_1(x)}{N(t, t_0, x)} \, dx.
\end{aligned}$$

Then, from (25) we obtain

$$\dot{w}(t, t_0, \omega) = (-\lambda_1 a_1 - \alpha + (1 - \epsilon) \gamma(\theta_t \omega)) w(t, t_0, \omega)$$

which leads to

$$w(t, t_0, \omega) \geq e^{-(\alpha + \lambda_1 a_1)(t - t_0) + (1 - \epsilon) \int_{t_0}^{t} \gamma(\theta_s \omega) ds} w(t_0), \ \forall t \geq \tau_0,$$

and (24) implies that $w(t, t_0, \omega) \to +\infty$ and $t \to +\infty$, which is a contradiction and therefore (22) holds true. $\square$

**Remark 4.1.** *In* (1)*-*(3)*, and let $\lambda_1$ be the first eigenvalue of $-\Delta$. If $a = 1$, and $m(\gamma)$ (see* (21)*) is independent of $\omega$, then it is possible to define the **basic reproductive number** as*

$$\mathcal{R}_0 := \frac{m(\gamma)}{\lambda_1 + d + b + c}.$$

*Thus, the disease is eradicated as long as $\mathcal{R}_0 < 1$, and the disease becames endemic $\mathcal{R}_0 > 1$, in the sense of* (22)*.*

**Example 4.2.** *Let $\lambda_1$ be the first eigenvalue of $-\Delta$, $a = 1$, $\gamma(\omega) = \gamma_0 + \Phi(\omega)$, with $\gamma_0 > 0$, and*

$$\lim_{t - t_0 \to +\infty} \frac{1}{t - t_0} \int_{t_0}^{t} \Phi(\theta_r \omega) \, dr = 0, \ \forall \omega \in \Omega.$$

*Then $m(\gamma) = \gamma_0$, and the basic reproductive number $\mathcal{R}_0 := \gamma_0 / (\lambda_1 + d + b + c)$. In [20] is given some examples of $\phi$ satisfying the above condition.*

## 5. Final Remarks

We close with a discussion of some open questions that emerge from our analysis and propose directions for future work.

1. All the results of this paper could be proven with weaker conditions for the recruitment rate $\Lambda$.

2. In the entire work, by dissipativness of the equation of total population, we could consider $SI$ in the model instead of $SI/N$.

3. If $\gamma$ is constant and $a_0 \leq a(x, \omega) \leq a_1$, then if $\gamma < \alpha_1 := \lambda_1 a_0 + d + b + c$ the disease is erradicated, and if $\gamma > \alpha_0 := \lambda_1 a_0 + d + b + c$, the infected persists in the sense of Theorem 4.2. Thus, if $\gamma \in [\alpha_0, \alpha_1]$ the behaviour is unkown.

4. We believe that the results of existence of attractors in Section 2, can be proven with environmental additive noise, as in [18].

## 6. Conclusion

In this work, we proposed a stochastic SIR model incorporating both random diffusion and random transmission effects to better capture the inherent uncertainties in disease dynamics. We provide insight into the long-term behavior of the system under stochastic influences, by studying the existence and properties of exponential attractors, random attractors, and special attracting sets. We examined how random transmission mechanisms lead to disease spread dynamics, drawing connections to non-autonomous bifurcation phenomena. Our approach departs from traditional energy-based methods and instead relies on the framework of mild solutions, allowing a less restrictive analysis. This offers a new approach to the study of random epidemic models and indicates directions for further research.

**Financial disclosure**

T. Caraballo and J. López-de-La-Cruz have been partially supported by the Spanish Ministerio de Ciencia e Innnovación, Agencia Estatal de Investigación (AEI) and Fondo Europeo de Desarrollo Regional (FEDER) under the project PID2024-??????-I00; A.N. Oliveira-Sousa has been partially supported by the São Paulo Research Foundation (FAPESP) under grant 2022/00176-0 and Coordenação de Aperfeiçoamento de Pessoal de Nível Superior - Brasil (CAPES) - Código de Financiamento 001 grants PROEX-9430931/D and 88887.912449/2023-00 (Print-CAPES); P. N. Seminario-Huertas has been partially supported by the Fondo Europeo de Desarrollo Regional (FEDER) under the project PID2021-122991NB-C21 and the Madrid Government (Comunidad de Madrid-Spain) under the Multiannual Agreement 2023-2026 with Universidad Politécnica de Madrid in the Line A, Emerging PhD researchers (grant agreement DOCTORES-EMERGENTES-24-UMDLRU-108-XVLWPC).